\documentclass{amsart}
\usepackage{amssymb,amsmath,amsfonts}
\usepackage{graphicx}
\usepackage{float}
\usepackage{multicol}

\newtheorem{theorem}{Theorem}[section]

\newtheorem{definition}{Definition}[section]

\numberwithin{equation}{section}
\begin{document}

\title[Minimal Degree Parameterizations of the Trefoil and Figure-Eight Knots]
{Minimal Degree Parameterizations of the trefoil and Figure-Eight Knots }
\author{Samantha Pezzimenti}
\address{Ramapo College of New Jersey, Mahwah, NJ 07430}
\email{spezzime@ramapo.edu}






\maketitle

\begin{abstract}
This paper determines the minimal degree sequence for two compact rational knots, namely the trefoil and figure-eight knots. We find explicit projections with the minimal degree sequence of each knot. This is done by modifying a non-compact rational minimal-degree parameterization of the trefoil and figure-eight knots to make it compact. This research was conducted at Ramapo College of New Jersey during the summer of 2011 with Dr. Donovan McFeron.
\end{abstract}

\pagestyle{myheadings} \markboth{Minimal Degree Parameterizations of the Trefoil and Figure-Eight Knots}{Minimal Degree Parameterizations of the Trefoil and Figure-Eight Knots}

\section{Introduction}
Determining the minimal degree sequence to parameterize a knot may be useful in classifying knots. We choose to work with compact rational knots because they can be contained in a finite region and can be expressed using rational functions. 

\smallskip
Along the lines of  \cite{DM2002}, we define,

\begin{definition}
Let $f(t), g(t), h(t)$ be three rational functions with degrees $\frac{p_{1}}{q_{1}}, \frac{p_{2}}{q_{2}}, \frac{p_{3}}{q_{3}}$ respectively, such that $q_{1} \leq q_{2} \leq q_{3}$, and if $q_{i} = q_{j}$ and $i<j$ then $p_{i} \leq p_{j}$. The degree sequence of $(f(t), g(t), h(t))$ is $\left(\frac{p_{1}}{q_{1}}, \frac{p_{2}}{q_{2}}, \frac{p_{3}}{q_{3}}\right)$. 
\end{definition}

\smallskip
Define the minimal degree sequence as,

\begin{definition}
{A degree sequence $\left(\frac{p_{1}}{q_{1}}, \frac{p_{2}}{q_{2}}, \frac{p_{3}}{q_{3}}\right)$ of a knot is minimal if for all other degree sequences $\left(\frac{l_{1}}{m_{1}}, \frac{l_{2}}{m_{2}}, \frac{l_{3}}{m_{3}}\right)$ of the knot, $\left(\frac{p_{1}}{q_{1}}, \frac{p_{2}}{q_{2}}, \frac{p_{3}}{q_{3}}\right) \leq  \left(\frac{l_{1}}{m_{1}}, \frac{l_{2}}{m_{2}}, \frac{l_{3}}{m_{3}}\right)$, which occurs when $(q_{1}, q_{2}, q_{3}, p_{1}, p_{2}, p_{3}) \leq (m_{1}, m_{2}, m_{3}, l_{1}, l_{2}, l_{3})$ using lexicographic ordering.}
\end{definition}

In  \cite{DM2002}, McFeron-Zuser established that the minimal degree sequence for the compact, rational trefoil and figure-eight knots is equal to or greater than $\left(\frac{2}{4}, \frac{2}{4}, \frac{2}{4}\right)$. They determine that the graph of the trefoil's double points must have four monotonic regions, which can only be represented by a degree $\frac{2}{4}$ function, in order to produce the correct over-under crossings. The same logic can be used for the figure-eight knot. 

This paper proves the theorems first conjectured in \cite{DM2002} by finding parameterizations for the trefoil and figure-eight knots with degree $\left(\frac{2}{4}, \frac{2}{4}, \frac{2}{4}\right)$.

\begin{theorem} The minimal degree sequence for the trefoil knot is $\left(\frac{2}{4}, \frac{2}{4}, \frac{2}{4}\right)$.
\end{theorem}

\begin{theorem} The minimal degree sequence for the figure-eight knot is $\left(\frac{2}{4}, \frac{2}{4}, \frac{2}{4}\right)$.
\end{theorem}

In section 2, we find a parameterization with the minimal degree sequence for the trefoil and explain the method in detail. We find a parameterization with the minimal degree sequence for the figure-eight knot in section 3 using a similar method. In section 4, we discuss the possibility of extending these methods for knots with five or more crossings. The calculations for this research were done primarily with Mathematica 8. 

\section{Minimal Degree Sequence for the Compact Rational Trefoil}
\subsection{$X$-$Y$ Projection of the trefoil}

As in  \cite{DM2002}, we let $f_{1}(t)$ and $g_{1}(t)$ be:
\begin{eqnarray}\label{ftrefoil}
f_{1}(t)&=&\frac{5t^{2}+2t+1}{1+t^2+t^3+t^4},
\end{eqnarray}
\begin{eqnarray}\label{gtrefoil}
g_{1}(t)&=&\frac{4t^{2}+1}{1+0.1t^2+t^4}.
\end{eqnarray}

\begin{figure}[H]
\begin{center}
\includegraphics[width=10cm]{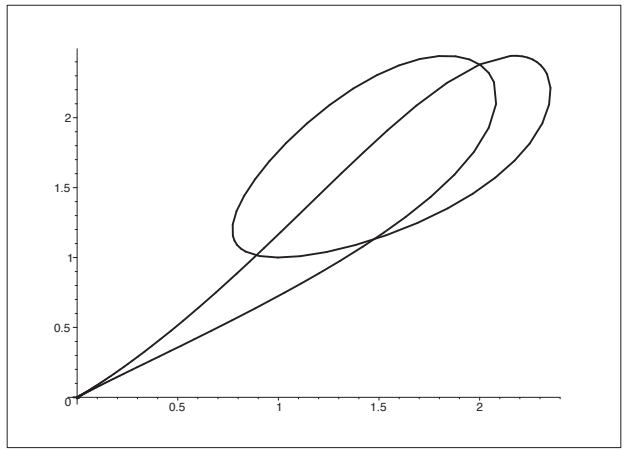}
\caption{Projection of the trefoil using the functions in McFeron-Zuser \cite{DM2002}
}
\end{center}
\end{figure}

\bigskip

This projection of the the trefoil has three double points such that for $(t_{i}<t_{i+1})$ we have:
\begin{eqnarray*}
f(t_{i})=f(t_{i+3}),
\\
g(t_{i})=g(t_{i+3}),
\\
i=1,2,3.
\end{eqnarray*}

The crossings of our $x$-$y$ projection occur approximately at:
\begin{eqnarray*}
t_{1} &=& -1.8461477068824201,
\\
t_{2} &=& -1,
\\
t_{3} &=& -0.06299423847617387,
\\
t_{4} &=& 0.1833157624425894,
\\
t_{5} &=& 1,
\\
t_{6} &=& 1.9583554137278836.
\end{eqnarray*}

We now attempt to find a $\frac{2}{4}$ function for $h_{1}(t)$ that produces alternating crossings.

\subsection{Finding $h_{1}(t)$ for the trefoil}

One way for $h_{1}(t)$ to produce alternating crossings on the trefoil is to require:
\begin{eqnarray*}
h(t_{1})>h(t_{4}),
\\
h(t_{2})<h(t_{5}),
\\
h(t_{3})>h(t_{6}).
\end{eqnarray*}

It is relatively straightforward to come up with a non-continuous degree $\frac{2}{4}$ function that satisfies the above condition, having four real roots in its denominator. From there, we will adjust the function to get rid of those real roots.

In order to create our temporary function, we will start by forming six linear factors with zeros occurring at points between each of the $t_{n}$ values. One example is $(t+2)(t+1.1)(t+.5)(t-.1)(t-.5)(t-1.5)$.

Since we want our function to be of degree $\frac{2}{4}$, we will pick two for the numerator and four for the denominator, giving us:

\begin{eqnarray}
h_{r_{1}}(t)=\frac{(t+2)(t+1.1)}{(t+.5)(t-.1)(t-.5)(t-1.5)}.
\end{eqnarray}

Expanding $h_{r_{!}}(t)$ gives us:

\begin{eqnarray*}
h_{r_{1}}(t)=\frac{t^2+3.1t+2.2}{t^4-1.6t^3-.1t^2+.4t-.0375}.
\end{eqnarray*}

Now we have a temporary function for $h$ that gives us that correct over-under crossings, but which has real zeros in the denominator. To get rid of the zeros, we will find the minimum of $t^4-1.6t^3-.1t^2+.4t-.0375$ and add a quantity greater than the minimum to the denominator. We find the minimum to be $-0.11608$ and will add $2$ to the denominator, making our new function:
\begin{eqnarray*}
h_{r_{2}}(t)=\frac{t^2+3.1t+2.2}{t^4-1.6t^3-.1t^2+.4t+1.9625}.
\end{eqnarray*}

Now we have no real roots in the denominator, but our over-under crossings have been altered. For $h_{r_{2}}$,
\begin{eqnarray*}
h(t_{1})<h(t_{4}),
\\
h(t_{2})<h(t_{5}),
\\
h(t_{3})<h(t_{6}).
\end{eqnarray*}

which is clearly not what we want. To correct this, we add a variable to one of our coefficients in $h_{r_{2}}(t)$ and plug in our $t_{n}$ values. We can then make three new functions:
\begin{eqnarray*}
 \gamma_{1}=h_{r_{2}}(t_{1})-h_{r_{2}}(t_{4}),
 \\
  \gamma_{2}=h_{r_{2}}(t_{2})-h_{r_{2}}(t_{5}),
 \\
  \gamma_{3}=h_{r_{2}}(t_{3})-h_{r_{2}}(t_{6}).
 \end{eqnarray*}

 We can find a point that satisfies our desired inequalities for $ \gamma_{1}$, $ \gamma_{2}$, and $ \gamma_{3}$. 

Through a series of similar adjustments, we can alter the other coefficients and end up with

\begin{eqnarray}\label{htrefoil}
h_{1}(t)=\frac{t^2+3.1t+3.2}{.981t^4+2.745t^3-.1t^2-1.6t+3.9685}.
\end{eqnarray}

We can now conclude that $\left(\frac{2}{4}, \frac{2}{4}, \frac{2}{4}\right)$ is the minimal degree sequence for the compact rational trefoil knot, thus proving Theorem 1.1.

\begin{figure}[H]
\begin{center}
\includegraphics[width=9cm]{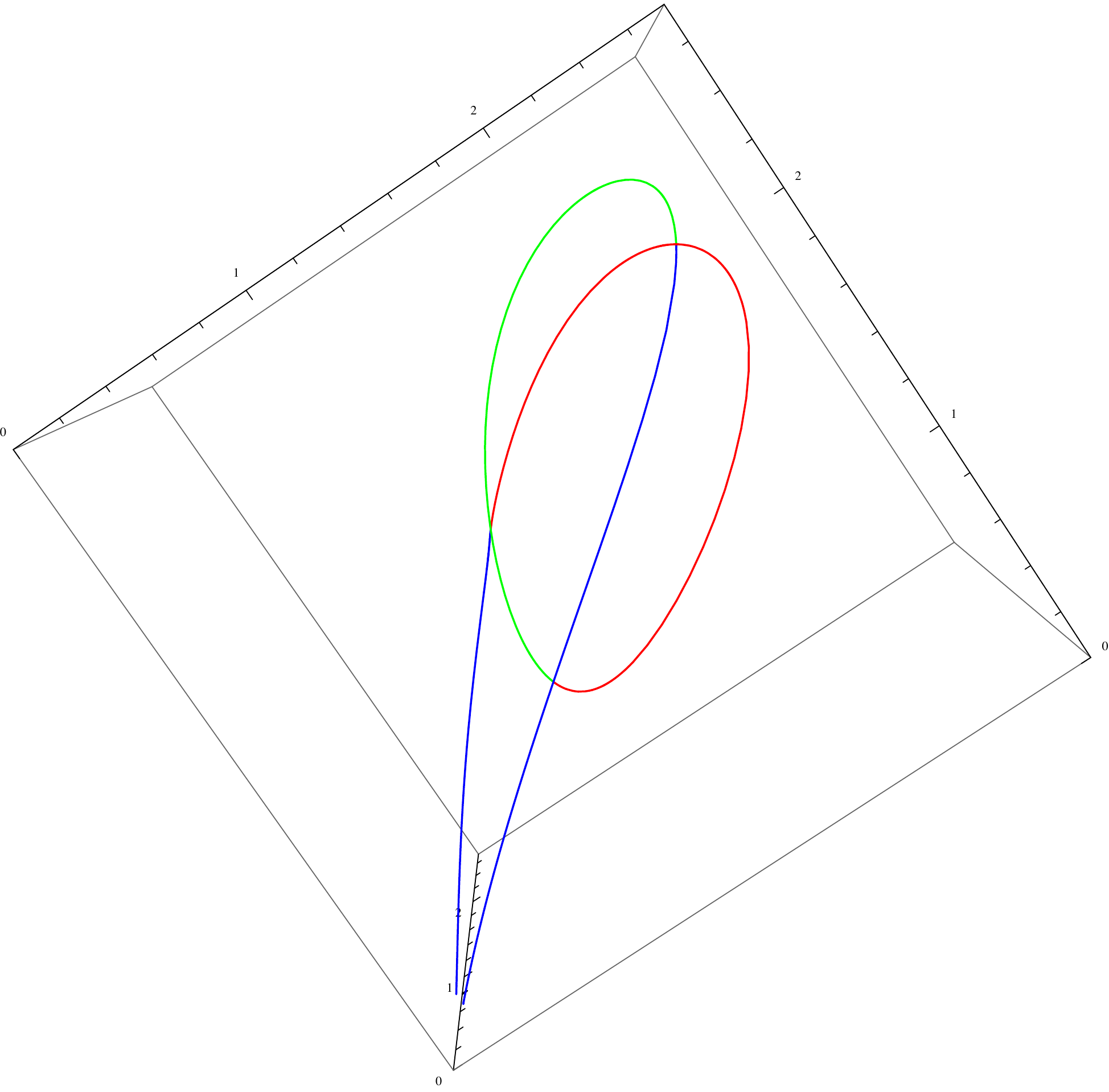}
\caption{The trefoil is tricolorable, as can be seen in our projection. A change in color signifies an under crossing \cite{CA2004}.}
\end{center}
\end{figure}

\section{Minimal Degree Sequence for the Compact Rational Figure-Eight Knot}

\subsection{Finding $h_{2}(t)$ for the figure-eight}

Using the algorithm defined in \cite{Algorithm} we find an $x$-$y$ projection of the figure-eight knot of degree $\left(\frac{6}{6}, \frac{6}{6}\right)$. We get
\begin{eqnarray*}
f_{\alpha}(t)&= &\frac{(t+7)(t+4)(t+1)(t-1)(t-4)(t-7)}{t^6+t^4+5000t^2+2392},
\\
g_{\alpha}(t)&=&\frac{(t+5)(t+3)(t-3)(t-5)(t-.972656)(t-7.027344)}{5096+t^2+23.514793t^4+t^6)}.
\end{eqnarray*}

\begin{figure}[H]
\begin{center}
\includegraphics[width=5cm]{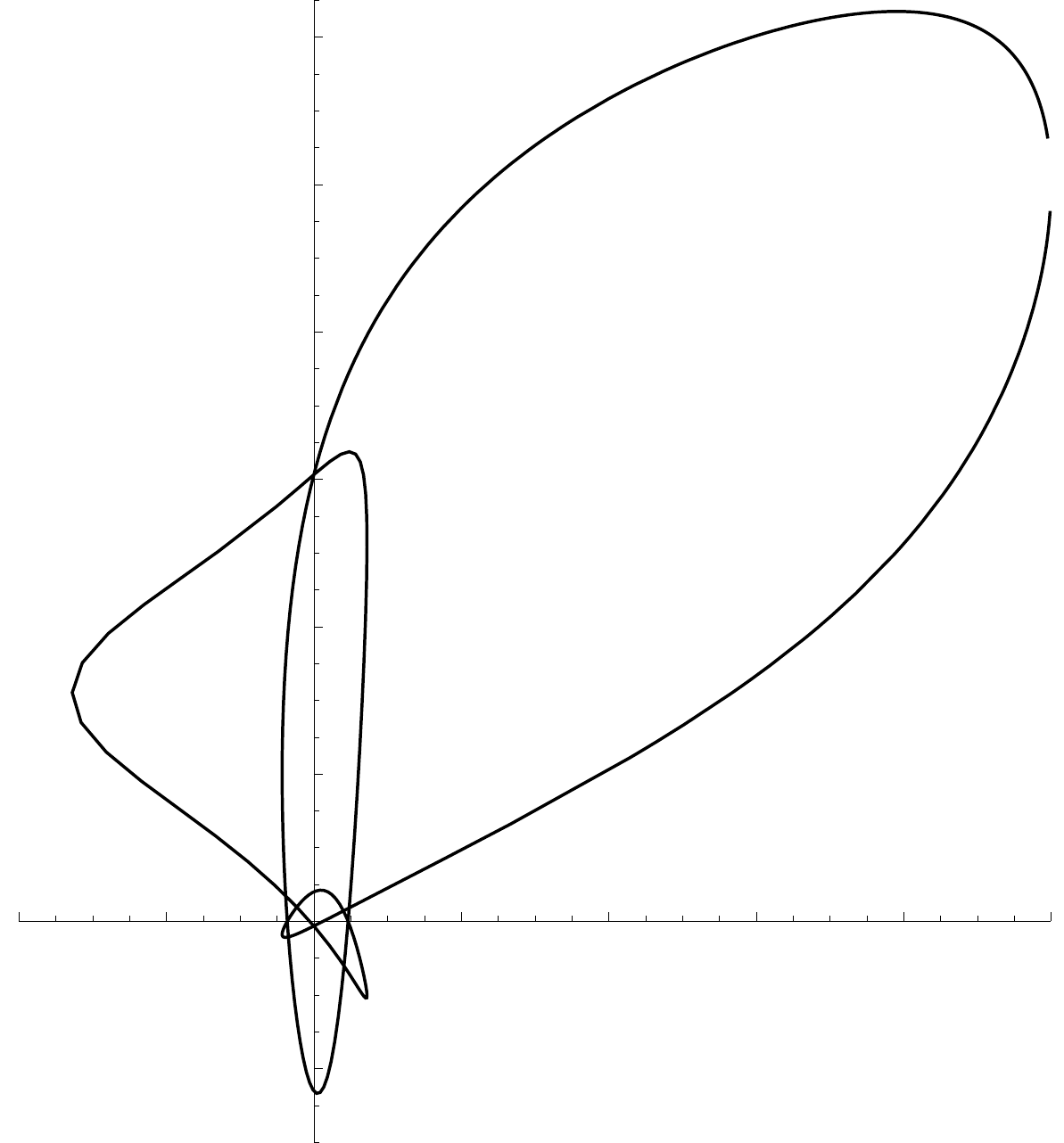}
\caption{$f_{\alpha}(t),g_{\alpha}(t)$}
\end{center}
\end{figure}

As with the trefoil knot, we find the points of intersection, that is, the points such that $f_{\alpha}(t_{i})=f_{\alpha}(t_{j})$ and $g_{\alpha}(t_{i})=g_{\alpha}(t_{j})$.  There are 20 double points which result in 10 crossings. We will call them $t_{1},t_{2},...,t_{20}$. Using these points, we determine which can be over crossings and which can be under crossings. We will use these sequences of over and under crossings to determine $h_{2}(t)$. There are many possible sequences that will result in a figure-eight knot. One example is


\begin{equation*} \begin{matrix}
h_{2}(t_{1})>h_{2}(t_{8}), & h_{2}(t_{6}) > h_{2}(t_{13}),
\\
h_{2}(t_{2})>h_{2}(t_{9}), & h_{2}(t_{7}) > h_{2}(t_{20}),
\\
h_{2}(t_{3}) < h_{2}(t_{16}), & h_{2}(t_{10}) < h_{2}(t_{15}), 
\\
h_{2}(t_{4}) > h_{2}(t_{17}), & h_{2}(t_{11}) < h_{2}(t_{18}), 
\\
h_{2}(t_{5}) > h_{2}(t_{12}), & h_{2}(t_{14}) > h_{2}(t_{19}).

\end{matrix} \end{equation*}

Using these inequalities, we find $h_{2}(t)$ with degree $\left(\frac{2}{4}, \frac{2}{4}\right)$  for our figure-eight Projection by applying the same method used to find $h_{1}(t)$ for the trefoil. We find $h_{2}(t)$ to be
\begin{eqnarray}\label{hfig8}
h_{2}(t)=\frac{-168.44 + 67.0899 t + 73.8617 t^2}{-0.484305 + t^4}.
\end{eqnarray}

\subsection{Finding $g_{2}(t)$ for the figure-eight}
Now, we will use $f_{\alpha}(t)$ and  $h_{2}(t)$ as our new $x$-$z$ projection and find the double points. We will use these points to determine the sequence of over and under crossings to use for $g_{2}(t)$.

\begin{figure}[H]
\begin{center}
\includegraphics[width=5cm]{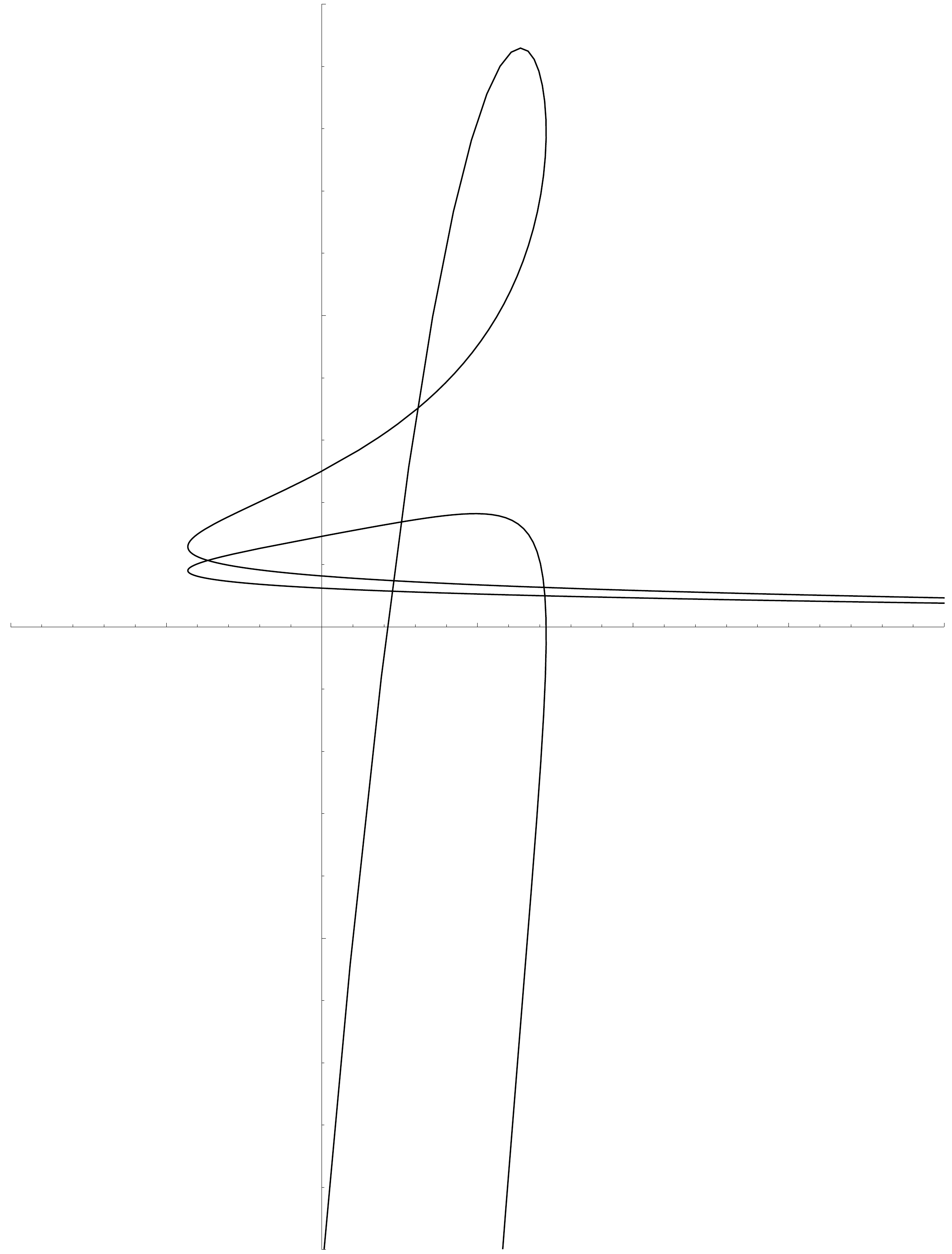}
\caption{$f_{\alpha}(t),h_{2}(t)$}
\end{center}
\end{figure}

We have 14 double points for the 7 crossings. Using the over under sequences we find 
\begin{eqnarray}\label{gfig8}
g_{2}(t)=\frac{-1591.53 - 455.993 t - 42.7391 t^2}{762.067 - 55.1785 t^2 + t^4}.
\end{eqnarray}

\subsection{Minimal degree sequence for the figure-eight}

Now, we can use $h_{2}(t), f_{2}(t)$ as our $x$-$y$ projection of the figure-eight knot and find $f_{2}(t)$ satisfying the correct over under sequence for the double points of our projection. We find, 
\begin{eqnarray}\label{ffig8}
f_{2}(t)=\frac{1.0516*10^{-12} + 4.72511*10^{-13} t + 4.9738*10^{-14} t^2}{785.103 + 
 29.5158 t - 24.8465 t^2 + 4.43299 t^3 + 0.960959 t^4}.
\end{eqnarray}

We now have a a degree $\left(\frac{2}{4}, \frac{2}{4}, \frac{2}{4}\right)$ that parameterizes the figure-eight knot, proving Theorem 1.2. We can conclude that the minimal degree sequence for the compact rational figure-eight knot is $\left(\frac{2}{4}, \frac{2}{4}, \frac{2}{4}\right)$.

\section{Knots with Greater than Four Crossings}
\subsection{Five Crossing Knots}

Exploring minimal degree sequences for knots with five or greater crossings is of interest for future research. We are currently exploring the minimal degree sequences for the two five crossing knots, which we will call $5_{1}$ and $5_{2}$. It can be seen that the correct over-under crossings for both the $5_{1}$ and $5_{2}$ knots can be achieved from a graph with four monotonic regions. 

For the $5_{1}$ knot, the $h(t)$ function must satisfy

\begin{eqnarray*}
h(t_{1})>h(t_{6}),
\\
h(t_{2})<h(t_{7}),
\\
h(t_{3})>h(t_{8}),
\\
h(t_{4})<h(t_{9}),
\\
h(t_{5})>h(t_{10}),
\end{eqnarray*}
and for the $5_{2}$ knot, the $h(t)$ function must satisfy

\begin{eqnarray*}
h(t_{1})>h(t_{6}),
\\
h(t_{2})<h(t_{5}),
\\
h(t_{3})>h(t_{8}),
\\
h(t_{4})<h(t_{9}),
\\
h(t_{7})>h(t_{10}).
\end{eqnarray*}

Figures 5 and 6 show that it is possible to achieve these over-under from a graph with four monotonic regions.

\begin{figure}[H]
\begin{center}
\includegraphics[width=9cm]{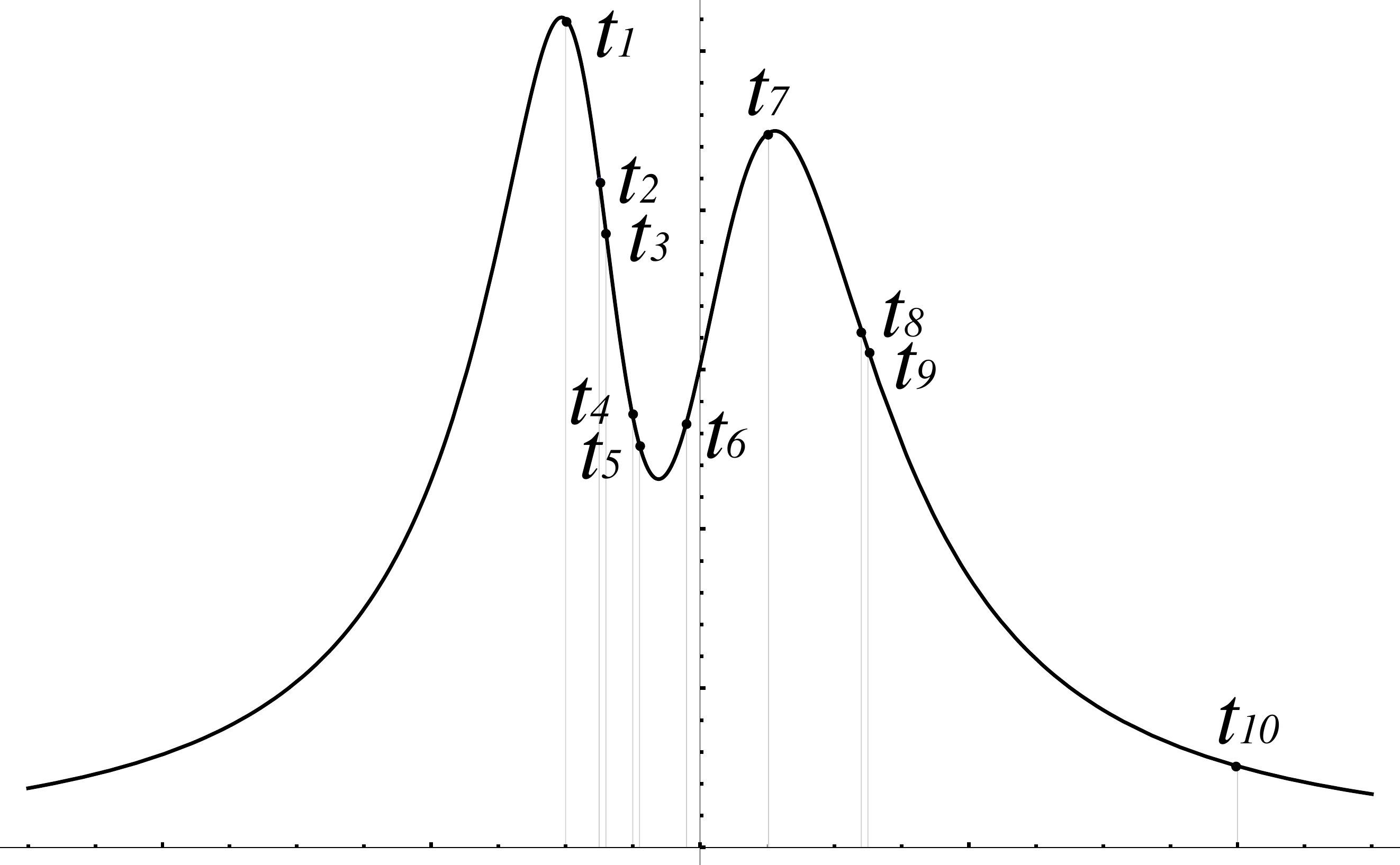}
\caption{Possible graph of the $5_{1}$ knot's crossings}
\end{center}
\end{figure}

\begin{figure}[H]
\begin{center}
\includegraphics[width=9cm]{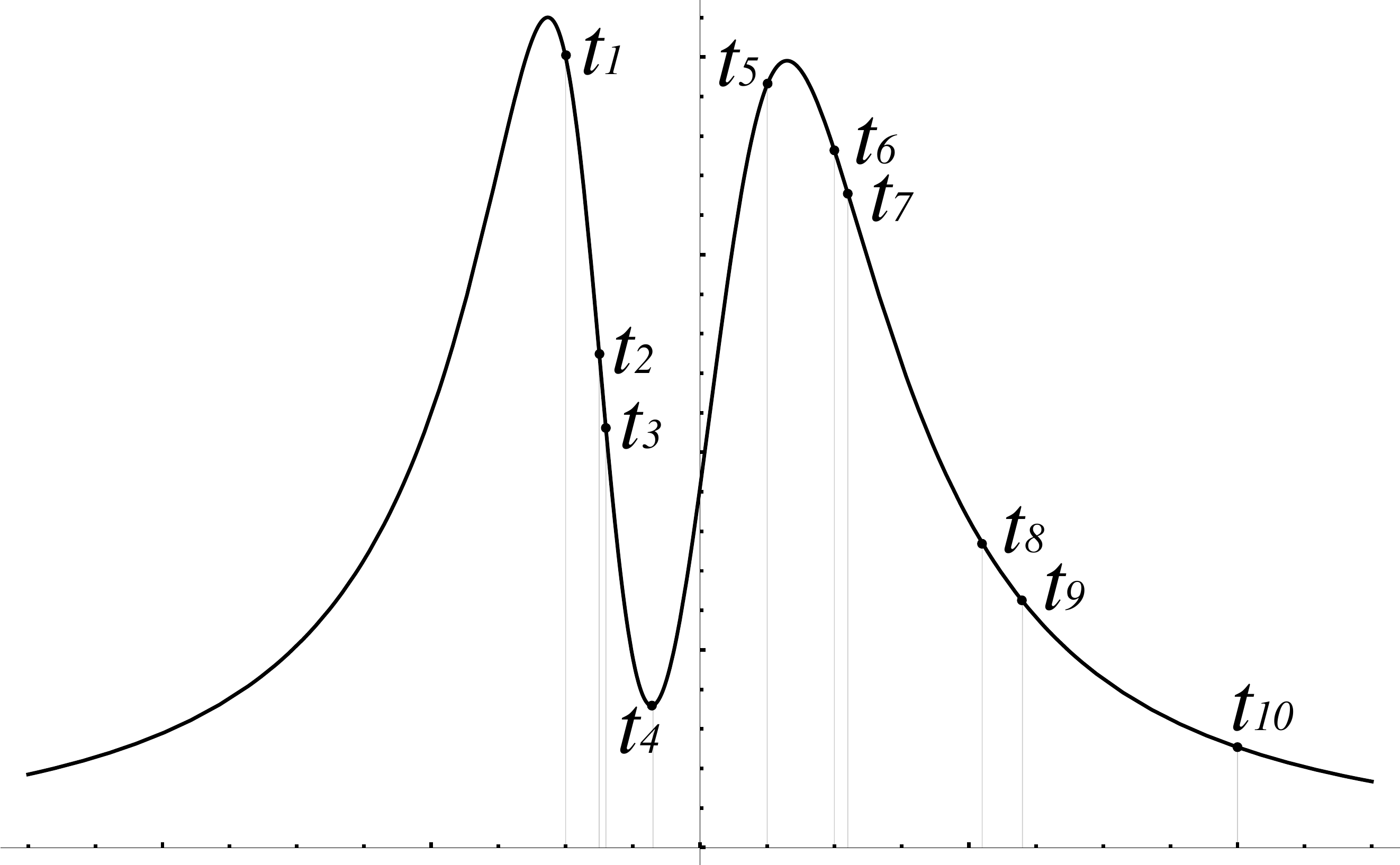}
\caption{Possible graph of the $5_{2}$ knot's crossings}
\end{center}
\end{figure}

As with the trefoil and figure-eight knots, this implies that the minimal degree sequence for a five crossing knot is greater than or equal to $\left(\frac{2}{4}, \frac{2}{4}, \frac{2}{4}\right)$.
Using a method similar to the one used to find the minimal degree functions for the figure-eight knot, we will first find minimal degree $x$-$y$ projections for each five crossing knot and then attempt to find $h(t)$ functions with degree $\frac{2}{4}$ that produce the correct crossings. 

%
%

We conclude that the minimal degree sequence for the compact rational  $5_{1}$ and $5_{2}$ knots is greater than or equal to $\left(\frac{2}{4}, \frac{2}{4}, \frac{2}{4}\right)$. 

\subsection{Knots with Greater than Five Crossings}
A long-term goal of this research is to classify knots based on their minimal degree sequence. Since we have only found minimal degree sequences equal to $\left(\frac{2}{4}, \frac{2}{4}, \frac{2}{4}\right)$, our first goal will be to find a knot whose minimal degree sequence is greater than $\left(\frac{2}{4}, \frac{2}{4}, \frac{2}{4}\right)$. 

\section{Conjectures}
\begin{itemize}
\item{The minimal degree sequence for the compact rational $5_{1}$ and $5_{2}$ knots is equal to $\left(\frac{2}{4}, \frac{2}{4}, \frac{2}{4}\right)$.}
\item{The minimal degree sequence for compact rational knots of six or more crossings is greater than $\left(\frac{2}{4}, \frac{2}{4}, \frac{2}{4}\right)$.}
\item{Knots of the same number of crossings will have the same minimal degree sequence.}
\end{itemize}



\begin{thebibliography}{10}
\bibitem{CA2004}
Colin C.~Adams.
\newblock The Knot Book.
\newblock American Mathematical Society. 2004.

\bibitem{DM2002}
Donovan ~McFeron, Alexandra ~Zuser.
\newblock On the Degrees of Rational Knots.
\newblock 2002.

\bibitem{Algorithm}
Donovan ~McFeron.
\newblock Algorithm for Constructing Compact Rational Knots.

\end{thebibliography}
\end{document}